\documentclass[11pt]{article}   

\usepackage{amsmath, amsthm, amsfonts, amssymb}
\usepackage{mathptmx}      
\usepackage{tikz}
\usetikzlibrary{shapes.multipart}
\definecolor{gold}{rgb}{1,.70,.0}   
\usepackage{float}                 

\newcommand{\QED}{\ $\square$}
\def\xxx{black}

\begin{document}

\title{\bf A note on unbounded Apollonian disk packings}

\author{Jerzy Kocik                  
\\ \small Department of Mathematics
\\ \small Southern Illinois University, Carbondale, IL62901
\\ \small jkocik{@}siu.edu  }

\date{\footnotesize (version 6 Jan 2019)}

\maketitle

\vspace{-.3in}

\begin{abstract}
\noindent
A construction and algebraic characterization of two unbounded Apollonian Disk packings 
in the plane and the half-plane are presented.  Both turn out to involve the golden ratio.
\end{abstract}

\section{Introduction}

We present two examples of unbounded Apollonian disk packings, 
one that fills a half-plane and one that fills the whole plane (see Figures \ref{fig:halfplane} and \ref{fig:plane}).  
Quite interestingly, both are related to the golden ratio.
We start with a brief review of Apollonian disk packings, define ``unbounded'', and fix notation and terminology.


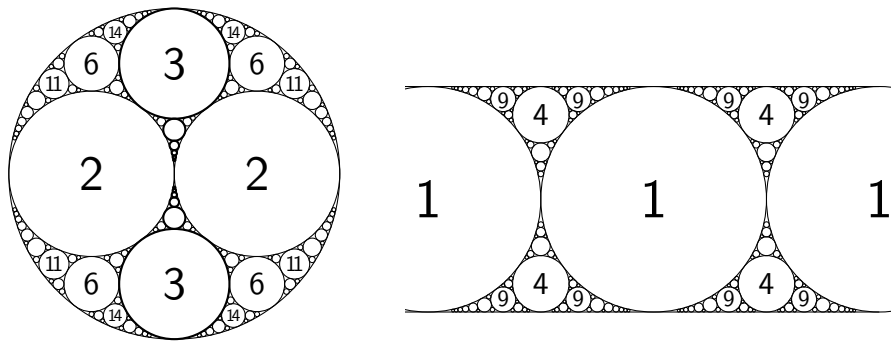
\begin{figure}[H]
\centering
\begin{tikzpicture}[scale=2.2]
\draw (0,0) circle (1);
\foreach \a/\b/\c   in {
1 / 0 / 2 
}
\draw (\a/\c,\b/\c) circle (1/\c)
          (-\a/\c,\b/\c) circle (1/\c);
\foreach \a/\b/\c in {
0 / 2 / 3 ,
0 /4 /15 ,
0 / 6 / 35 , 
0 / 8/ 63,
0 /10 / 99,
0 / 12 / 143
}
\draw[thick] (\a/\c,\b/\c) circle (1/\c)
          (\a/\c,-\b/\c) circle (1/\c) ;

\foreach \a/\b/\c/\d in {
3 / 4 /6, 	8 / 6 / 11,	5 / 12/ 14,	15/ 8 / 18,	8 / 12 / 23,	7 / 24 / 26,
24/	10/	27, 	21/	20/	30, 	16/	30/	35, 	3/	12/	38, 	35/	12/	38, 	24/	20/	39, 	9/	40/	42,
16/	36/	47, 	15/	24/	50, 	 48/	14/	51, 	45/	28/	54, 	24/	30/	59, 	40/	42/	59, 	11/	60/	62,
21/	36/	62, 	48/	28/	63, 	33/	56/	66, 	63/	16/	66, 	8/	24/	71, 	55/	48/	74, 	24/	70/	75,
48/	42/	83, 	80/	18/	83, 	13/	84/	86, 	77/	36/	86, 	24/	76/	87, 	24/	40/	87, 	39/	80/	90
}
\draw (\a/\c,\b/\c) circle (1/\c)       (-\a/\c,\b/\c) circle (1/\c)
          (\a/\c,-\b/\c) circle (1/\c)       (-\a/\c,-\b/\c) circle (1/\c) 
;
\node at (-1/2,0) [scale=1.7, color=\xxx] {\sf 2};
\node at (1/2,0) [scale=1.7, color=\xxx] {\sf 2};
\node at (0,2/3) [scale=1.6, color=\xxx] {\sf 3};
\node at (0,-2/3) [scale=1.6, color=\xxx] {\sf 3};
\node at (1/2,2/3) [scale=1.1, color=\xxx] {\sf 6};
\node at (-1/2,2/3) [scale=1.1, color=\xxx] {\sf 6};
\node at (1/2,-2/3) [scale=1.1, color=\xxx] {\sf 6};
\node at (-1/2,-2/3) [scale=1.1, color=\xxx] {\sf 6};
\node at (8/11,6/11) [scale=.77, color=\xxx] {\sf 1$\!$1};
\node at (-8/11,6/11) [scale=.77, color=\xxx] {\sf 1\!1};
\node at (8/11,-6/11) [scale=.77, color=\xxx] {\sf 1\!1};
\node at (-8/11,-6/11) [scale=.77, color=\xxx] {\sf 1\!1};
\node at (5/14,6/7) [scale=.6, color=\xxx] {\sf 1\!4};
\node at (-5/14,6/7) [scale=.6, color=\xxx] {\sf 1\!4};
\node at (5/14,-6/7) [scale=.6, color=\xxx] {\sf 1\!4};
\node at (-5/14,-6/7) [scale=.6, color=\xxx] {\sf 1\!4};
\end{tikzpicture}
\qquad
%
%
\begin{tikzpicture}[scale=1.5, rotate=90, shift={(0,2cm)}]  
\clip (-1.25,-2.1) rectangle (1.1,2.2);
\draw (1,-2) -- (1,3);
\draw (-1,-2) -- (-1,3);
\draw (0,0) circle (1);
\draw (0,2) circle (1);
\draw (0,-2) circle (1);

\foreach \a/\b/\c in {
3/4/4,  5/12/12,  7/24/24, 9/40/40  
}
\draw (\a/\c,\b/\c) circle (1/\c)    (\a/\c,-\b/\c) circle (1/\c)
          (-\a/\c,\b/\c) circle (1/\c)    (-\a/\c,-\b/\c) circle (1/\c)   ;

\foreach \a/\b/\c in {
8/    6/   9, 	
15/   8/   16 , 	
24/  20/  25, 	
24/  10/  25, 	
21/  20/  28,
16/	30/	33,    
35/	12/	36,
48/	42/	49,
48/	28/	49,
48/	14/	49,
45/	28/	52,
40/	42/	57,
33/	56/	64,
63/	48/	64,
63/	16/	64,
55/	48/	72,
24/	70/	73,
69/	60/	76,
80/	72/	81,
64/	60/	81,
80/	36/	81,
80/	18/	81
}
\draw (\a/\c, \b/\c) circle (1/\c)          (-\a/\c, \b/\c) circle (1/\c)
          (\a/\c,-\b/\c) circle (1/\c)         (-\a/\c,-\b/\c) circle (1/\c)
          (\a/\c,2-\b/\c) circle (1/\c)       (-\a/\c,2-\b/\c) circle (1/\c)
          (\a/\c,2+\b/\c) circle (1/\c)       (-\a/\c,2+\b/\c) circle (1/\c)
          (\a/\c,-2+\b/\c) circle (1/\c)       (-\a/\c,-2+\b/\c) circle (1/\c)
;
\node at (0,0) [scale=1.9, color=\xxx] {\sf 1};
\node at (0,-2) [scale=1.9, color=\xxx] {\sf 1};
\node at (0,2) [scale=1.9, color=\xxx] {\sf 1};
\node at (-3/4,1) [scale=1.1, color=\xxx] {\sf 4};
\node at (-3/4,-1) [scale=1.1, color=\xxx] {\sf 4};
\node at (3/4,1) [scale=1.1, color=\xxx] {\sf 4};
\node at (3/4,-1) [scale=1.1, color=\xxx] {\sf 4};
\node at (7/8,8/6) [scale=.8, color=\xxx] {\sf 9};
\node at (7/8,-8/6) [scale=.8, color=\xxx] {\sf 9};
\node at (7/8,4/6) [scale=.8, color=\xxx] {\sf 9};
\node at (7/8,-4/6) [scale=.8, color=\xxx] {\sf 9};
\node at (-7/8,8/6) [scale=.8, color=\xxx] {\sf 9};
\node at (-7/8,-8/6) [scale=.8, color=\xxx] {\sf 9};
\node at (-7/8,4/6) [scale=.8, color=\xxx] {\sf 9};
\node at (-7/8,-4/6) [scale=.8, color=\xxx] {\sf 9};
\end{tikzpicture}
\caption{Apollonian Window (left) and Apollonian Belt (right).}
\label{fig:Apollo}
\end{figure}

The Apollonian disk packing
is a fractal arrangement of disks such that
any of its three mutually tangent disks determine it  by recursivly inscribing new disks 
in the curvy-triangular spaces that emerge between the disks.
In such a context, the three initial disks are called a {\bf seed} of the packing.
Figure~\ref{fig:Apollo} shows for two most popular examples:
the ``Apollonian Window'' and the ``Apollonian Belt''.
The numbers inside the circles represent their curvatures (reciprocals of the radii).

Note that the curvatures are integers; such arrangements are called {\bf integral} 
Apollonian disk packings; they are  classified and their properties are still studied \cite {jk,LMW,N}. 
Since much effort has been invested in the study of integral Apollonian packings because 
of their connections with number theory, 
one could possibly  get the impression that they represent all types of configurations.

Na{\"\i}vely, the 
Apollonian Window seems bounded and the Apollonian Belt does not.
But this is not so.
Actually both arrangements cover the whole plane.
This is because the most external circle in the Apollonian Window is the  boundary of an infinite
external disk; such a disk is considered to have  negative radius (and curvature).
Similarly, the two lines in the Apollonian belt are actually 
circles of zero curvature and they bound disks that are tantamount to two half-planes.
They are tangent at infinity.
And again, the whole plane is covered.  Hence the definition:

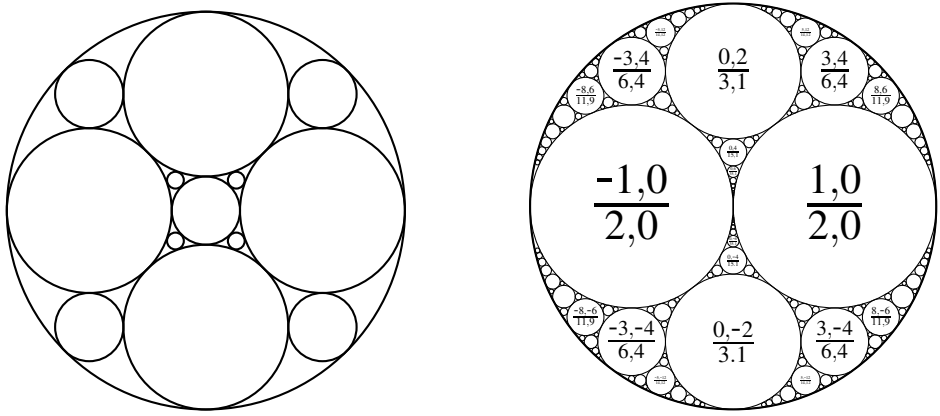
\begin{figure}[t]
\centering
\begin{tikzpicture}[scale=1.55]
\draw [thick] (0,0) circle (1.70711);
\draw [thick] (1,0) circle (.70711);
\draw [thick] (0,1) circle (.70711);
\draw [thick] (-1,0) circle (.70711);
\draw [thick] (0,-1) circle (.70711);
\draw [thick] (0,0) circle (.29289);
\draw [thick] (.261,.261) circle (.072);
\draw [thick] (-.261,.261) circle (.072);
\draw [thick] (.261,-.261) circle (.072);
\draw [thick] (-.261,-.261) circle (.072);
\draw [thick] (1,1) circle (.29289);
\draw [thick] (1,-1) circle (.29289);
\draw [thick] (-1,1) circle (.29289);
\draw [thick] (-1,-1) circle (.29289);
\end{tikzpicture} \qquad\qquad
\begin{tikzpicture}[scale=2.7]
\draw (0,0) circle (1);
\draw [thick] (0,0) circle (1);
\foreach \a/\b/\c/\d/\s   in {
1 / 0 / 2 / 0 / 2
}
\draw (\a/\c,\b/\c) circle (1/\c)
          (-\a/\c,\b/\c) circle (1/\c)
node  [scale=\s]  at (\a/\c,\b/\c)   {$\frac{\a,\b}{\c,\d}$}
node  [scale=\s]  at (-\a/\c,\b/\c) {$\frac{\hbox{-}\a,\b}{\c,\d}$};

\foreach \a/\b/\c/\d/\s in {
0 / 2 / 3  /1 /1.1,
0 /4 /15 /1  /.3,
0 / 6 / 35 / 1/ .2 
}
\draw (\a/\c,\b/\c) circle (1/\c)
          (\a/\c,-\b/\c) circle (1/\c)
        node [scale= \s]  at (\a/\c,\b/\c)  {$\frac{\a,\b}{\c,\d}$}
        node [scale= \s]  at (\a/\c,-\b/\c) {$\frac{\a,\hbox{-}\b}{\c.\d}$};

\foreach \a/\b/\c/\d in {
0 / 8/ 63 / 1,
0 /10 / 99 / 1
}
\draw (\a/\c,\b/\c) circle (1/\c)  
           (\a/\c,-\b/\c) circle (1/\c) 
;

\foreach \a/\b/\c/\d in {
0 / 12 / 143 / 1,
0 / 14 / 195 / 1,
0 / 16 / 255 / 1,
0 / 18 / 323 / 1
}
\draw (\a/\c,\b/\c) circle (1/\c)  
           (\a/\c,-\b/\c) circle (1/\c);

\foreach \a/\b/\c/\d/\s in {
3 / 4 / 6/ 4 / 1,
8 / 6 / 11/ 9 / .5,
5 / 12/ 14/ 12 / .2
}
\draw (\a/\c,\b/\c) circle (1/\c)       (-\a/\c,\b/\c) circle (1/\c)
          (\a/\c,-\b/\c) circle (1/\c)       (-\a/\c,-\b/\c) circle (1/\c)
node [scale=\s]  at (\a/\c,\b/\c)    {$\frac{\a, \b}{\c,\d}$}
node [scale=\s]  at (-\a/\c,\b/\c)   {$\frac{\hbox{-}\a,\b}{\c,\d}$}
node [scale=\s]  at (\a/\c,-\b/\c)   {$\frac{\a,\hbox{-}\b}{\c,\d}$}
node [scale=\s]  at (-\a/\c,-\b/\c)  {$\frac{\hbox{-}\a,\hbox{-}\b}{\c,\d}$};

\foreach \a /  \b / \c / \s  in {
15/ 8 / 18 / 16,
8 / 12 / 23 / 9,
7 / 24 / 26 / 24
}
\draw (\a/\c,\b/\c) circle (1/\c)          (-\a/\c,\b/\c) circle (1/\c)
          (\a/\c,-\b/\c) circle (1/\c)         (-\a/\c,-\b/\c) circle (1/\c)
;


\foreach \a /  \b / \c / \d  in {
24/	10/	27/	25,
21/	20/	30/	28,
16/	30/	35/	33,
3/	12/	38/	4,
35/	12/	38/	36,
24/	20/	39/	25,
9/	40/	42/	40,
16/	36/	47/	33,
15/	24/	50/	16,
48/	14/	51/	49,
45/	28/	54/	52,
24/	30/	59/	25,
40/	42/	59/	57,
11/	60/	62/	60,
21/	36/	62/	28,
48/	28/	63/	49,
33/	56/	66/	64,
63/	16/	66/	64,
8/	24/	71/	9,
55/	48/	74/	72,
24/	70/	75/	73,
48/	42/	83/	49,
80/	18/	83/	81,
13/	84/	86/	84,
77/	36/	86/	84,
24/	76/	87/	73,
24/	40/	87/	25,
39/	80/	90/	88,
64/	60/	95/	81,
80/	36/	95/	81,
33/	72/	98/	64,
63/	48/	98/	64,
65/	72/	98/	96
}
\draw (\a/\c, \b/\c)  circle (1/\c)         (-\a/\c,  \b/\c)  circle  (1/\c)
          (\a/\c, -\b/\c) circle (1/\c)         (-\a/\c, -\b/\c)  circle  (1/\c);
\end{tikzpicture}
\caption{Left: This is not an Apollonian packing; \ Right: Apollonian window with symbols.}
\label{fig:two}
\end{figure}


~\\
{\bf Definition:}   An Apollonian packing is called {\bf unbounded} if 
for any $x\in \mathbb R$,  it contains a disk 
with radius $r$ satisfying $x<r<\infty$. 
\\

Note that packings like the ones in Figure~\ref{fig:Apollo} are not unbounded.
The maximal radius in the first case is 1/2,
and the maximal (finite) radius in the second case is 1.

~

Let us recall also a few facts concerning  disks in the coordinated Euclidean plane $E\cong \mathbb R^2$.
For a disk with  center at $(x,y)$ and radius $r$, we introduce
``{\bf reduced coordinates}'' $\dot x=x/r$ and $\dot y= y/r$.  The curvature is denoted by $\beta = 1/r$.
Disks will be coded by {\bf symbols}, formal fractions that are related to standard coordinates as follows:
%
$$
\hbox{Circle}\left( (x,y),\, r\right) \ \Rightarrow 
\frac{\dot x,\ \dot y}{\beta,\ \gamma} \qquad \hbox{so that}\qquad 
(x,y) = \left(\frac{\dot x}{\beta}, \ \frac{\dot y}{\beta} \right), \  \ \ r=1/\beta 
$$
The fourth term, $\gamma$, is the``co-curvature'', 
which represents the curvature of the image of the disk under inversion in the unit circle. 
(It will not play a crucial  role in these notes.)
The symbol is a compact representation of a vector in the 4-dimensional linear space $M\cong\mathbb R^4$:
$$
\frac{\dot x,\ \dot y}{\beta,\ \gamma}
 \quad \equiv \quad 
\begin{bmatrix}
\dot x\\  \dot y\\ \beta\\ \gamma
\end{bmatrix}
$$
The space $M$ is equipped with a (billinear) inner product defined
for two disks $D_1$ and $D_2$ as the real number
\begin{equation}
\label{eq:inner}
\langle D_1, D_2\rangle  \ = \
-\dot x_1 \dot x_2 - \dot y_1 \dot y_2+ \frac {1}{2}\left (\beta_1\gamma_2 + \gamma_1\beta_2\right) 
\end{equation}
In particular, two disks are tangent (externally) if 
\begin{equation}
\label{eq:tangent}
\langle D_1, D_2\rangle = 1
\qquad\hbox{or}\qquad
-\dot x_1 \dot x_2 - \dot y_1 \dot y_2+ \frac {1}{2}\left (\beta_1\gamma_2 + \gamma_1\beta_2\right) = 1
\end{equation}
The norm of any disk is 
\begin{equation}
\label{eq:norm}
\|D\|^2 \ \equiv \ \langle D, D\rangle = -1
\qquad\hbox{or}\qquad
-\dot x^2 - \dot y^2+ \beta\gamma = -1
\end{equation}
(Note that this is sufficient to determine the value of the co-curvature $\gamma$ without much ado.)

~

Four circles are in {\bf Descartes configuration} if they are pairwise tangent to each other. 
In such a case their curvatures satisfy the {\bf Descartes formula}:
\begin{equation}
\label{eq:Descartes}
2\,(a^2 + b^2 + c^2 + d^2) = (a+b+c+d)^2 \,.
\end{equation}
This equation is part of an extended Descartes formula that has the following matrix form:
\begin{equation}
\label{eq:Descartes2}
MFM^T=G \,,
\end{equation}
where
$$
M=\begin{bmatrix}
\dot x_1 & \dot x_2 & \dot x_3 &\dot x_4\\ 
\dot y_1 & \dot y_2 & \dot y_3 &\dot y_4\\ 
\beta_1 & \beta_2 & \beta_3 &\beta_4\\ 
\gamma_1 & \gamma_2 & \gamma_3 &\gamma_4\\ 
\end{bmatrix}\,,
\quad
F=\left[\begin{array}{rrrr}
-1 & 1 & 1&1 \\ 
1 & -1 & 1&1 \\
 1 & 1 & -1&1 \\ 
1 & 1 & 1&-1  
\end{array}
\right]\,,
\quad
G=\begin{bmatrix}
-4 & 0 & 0&0 \\ 
0 & -4 & 0&0 \\
 0 & 0 & 0&2 \\ 
0 & 0 & 2&0  
\end{bmatrix}\,.
$$
For a proof see \cite{LMW,jk1}. ( The proof in \cite{jk1} is based on the concepts presented above.)
Three mutually tangent circles may be complemented by a fourth circles to form a Descartes configuration in two ways.
One may solve the quadratic equations \eqref{eq:Descartes},
or the more general  \eqref{eq:Descartes2},
to obtain these solutions:
$$
\begin{array}{rl}
\dot x_4 =& \dot x_1 + \dot x_2 + \dot x_3 
                       \  \pm \ \sqrt{ \dot x_1\dot x_2 + \dot x_2 \dot x_3 + \dot x_3 \dot x_1 +1}\\[7pt]
\beta_4 =& \beta_1 + \beta_2 + \beta_3 
                      \  \pm \ \sqrt{ \beta_1\beta_2 + \beta_2 \beta_3 + \beta_3 \beta_1}
\end{array}
$$
The expressions for $\dot y$ and $\gamma$
are analogous to these  for $\dot x$ and $\beta$, respectively. 
The ``$\pm$'' captures both solutions to the problem --- two possible constructions of the fourth disk.

~

Finally, a few comments on the golden means.
We distinguish between the {\bf golden ratio} $\varphi$ and the {\bf golden cut} $\tau$:
\begin{equation}
\label{eq:phi}
\varphi = \frac{1+\sqrt{5}}{2} = 1.618...
\qquad \hbox{and}\qquad
\tau = \frac{-1+\sqrt{5}}{2} = 0.618... \ ,
\end{equation}
Recall that  $\varphi = \tau +1$,  $\varphi  \tau  =1 $ and $\tau^k = \varphi^{-k}$.
The Fibonacci numbers are understood as a bilateral sequence $(F_i)$ defined by $F_0=0$, $F_1=1$,
and the recurrence $F_{i+1} = F_i+F_{i+1}$.  The labeling goes as follows:
$$
\begin{array}{cccccccccccccccc}
\ldots &F_{-6}&F_{-5} &F_{-4} &F_{-3} &F_{-2} &F_{-1} 
                    &F_{0} & F_{1} &F_{2} &F_{3} &F_{4} &F_{5} &F_{6} &F_{7} &\ldots \\
\dots&-8&5&-3&2&-1&1&0&1&1&2&3&5&8&13&\ldots
\end{array}
$$
One notices that $F_{-n} = (-1)^n F_n$.
The Fibonacci numbers are related to the golden means by a formula that is given here in two forms
\begin{equation}
\label{eq:Fphi}
\varphi^n \ = \ F_n\varphi +F_{n-1}
\qquad\hbox{and}\qquad 
\tau^n \ = \ (-1)^n \left(F_{n+1} - F_{n+1}\varphi\right)
\end{equation}
(The second formula results from the first under replacement $n\mapsto -n$ and $\varphi=1+\tau$.)

\section{Half-plane packing}

\begin{figure}[t]
\begin{center}
\begin{tikzpicture}[scale=2]

\clip (-3.5, -.05) rectangle  (3.5,3.5);

\draw (-4,0) -- (4,0);


\foreach \a/\b   in {
8.472135955/	8.972135955,
-2.618033989/	3.427050983,
1.618033989/	1.309016994,
0/	0.5,
0.618033989/	0.190983006,
0.381966011/	0.072949017,
0.472135955/	0.027864045,
0.437694101/	0.010643118,
0.450849719/	0.004065309
}
\draw[fill=gold!10]  (\a,\b) circle (\b)
;
\node at (2.8,2.8) [align=left] {\Large $d = \varphi^6$}; 
\node at (-2.618033989+1,	3.427050983-2) [align=left] {\Large$d=\varphi^4$ }; 
\node at (1.618033989, 1.309016994)  [align=left] {\Large$d=\varphi^2$ }; 
\node at (0,	0.5) [align=left] {\large$d=1$ }; 
\node at (0.618033989, 0.190983006)[align=left, scale=.9] {$d\!\!=\!\!\tau^{2}$}; 

%
\end{tikzpicture}
\end{center}
\caption{Upper half-plane Apollonian disk packing.  Diameters are denoted as $d$.}
\label{fig:halfplane}
\end{figure}
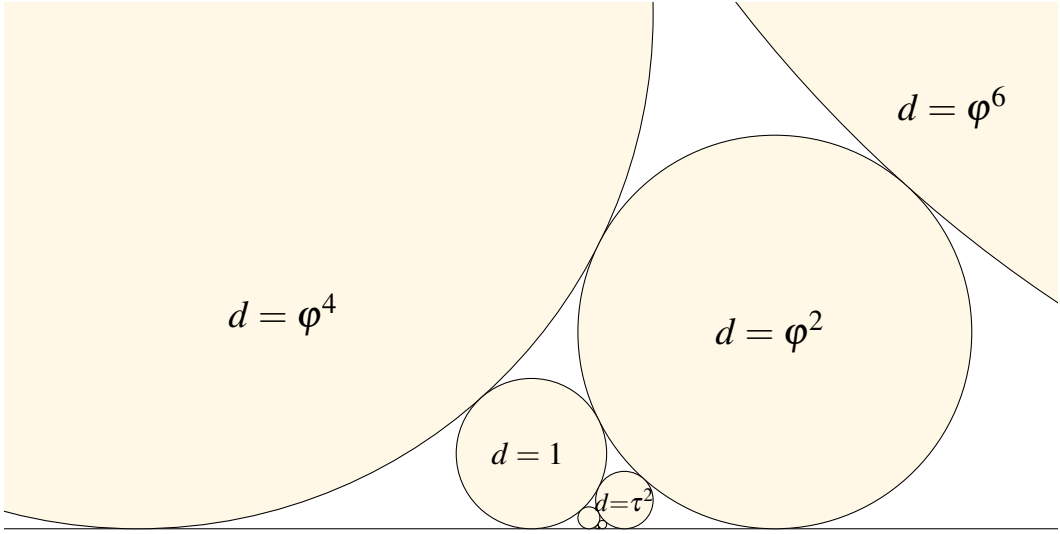

An example of an unbounded Apollonian disk packing filling the upper half-plane is shown in Figure~\ref{fig:halfplane}. 
Only the central zigzag chain of ``golden disks'' is displayed.
The zigzag consists of a sequence of disks of diameters equal to the even powers of the golden ratio $\varphi$.  
Any three consecutive disks are mutually tangent and all are tangent to the horizontal line identified with the $x$-axis.  
The following proposition assures existence of such an arrangement.

~\\ \noindent
{\bf Proposition 1:}  The sequence $D_n$ of disks with symbols given by   
\begin{equation}
\label{eq:HalfGold}
\mathcal D_n  =  \displaystyle\frac{2F_n\varphi^n,\; 1}{~~~2\varphi^{2n},\;2F_n^2}\,,  \qquad n\in \mathbb Z
\end{equation}
provides an infinite skeleton for an Apollonian arrangement filling the upper half-plane
with disks of arbitrarily large radius.
Every three consecutive disks and the horizontal line $y=0$ form a Descartes configuration.
In particular, the following three disks
$$
\left(~
 \dfrac{0,\; -1}{0,\; ~~0}\,,  \quad 
  \dfrac{0,\; 1}{2,\;0}\,,  \quad 
   \dfrac{2\varphi,\; 1}{2\varphi^{2},\;2}\
~\right)
$$
may serve as a seed of this arrangement.

~\\
{\bf Proof:} Applying the extended Descartes formula \eqref{eq:Descartes2}
 for the four consecutive entries of the above terms would suffice.
Or one may simply verify that each of the disk in the chain is tangent to the axis line
and that any two consecutive, and step-2 consecutive, disks  are mutually tangent:
$$
\langle  D_n,\,  D_{n+1} \rangle = 1 
\qquad\hbox{and} \qquad
\langle D_{n},\, D_{n+2} \rangle = 1
$$
The first condition, using the explicit expression for the inner product \eqref{eq:inner}, may be written as
$$
-2F_nF_{n+1}\varphi^{2n+1} + F_{n+1}^2 \varphi^{2n} + F^2_n\varphi^{2n+2} = 1
$$
The left hand side may be reduced as follows:
$$
\varphi^{2n} \left( F_{n+1} - \varphi F_n\right)^2
\ = \ 
\varphi^{2n} \left( \pm \tau^n\right)^2 
\ = \ 
1\,,
$$
as required.  The second property follows similarly.
\QED
 
%

~\\
The disk {\bf diameters} ($d=2r$) and their horizontal {\bf positions} $x$
may be read off from \eqref{eq:HalfGold}:
$$
\begin{array}{llcrrrrrrrrl}
d_n = \varphi^{-2n}   
&:
&\quad \ldots &\varphi^{6}, 
&  \varphi^{4}, &  \varphi^{2}, &  1, &  \tau^2, &  \tau^4, &  \tau^6, &  \tau^8,\ldots 
&\hbox{\small ~~towards smaller disks}   \\
x_n = F_n\varphi^{-n} = F_n\tau^{n}
&:
&\quad \ldots &2 \varphi^{3},
&  -\varphi^{2}, &  \varphi, &  0, &  \tau, &  \tau^2, &  2\tau^3, & 3\tau^4,\ldots 
&\hbox{\small ~~approaching } 5^{-1/2}
\end{array}
$$
(Index $n$ runs from $-\infty$ to $+\infty$, from larger disks to smaller.)  
The unboundedness of the arrangement is thus evident.
As we follow the sequence towards the smaller circles, the points of tangency to the axis line converge to the limit point:
$$
x_{\infty} \ = \ \lim_{n\to\infty} x_n
 \ = \ \lim_{n\to\,\infty} \frac{F_n}{\varphi^n}
  \ = \ \frac{1}{\sqrt{5}} \,.
$$

\begin{figure}[t]
\begin{center}
\begin{tikzpicture}[scale=2.1]

\clip (-3.5, -.79) rectangle  (3.2,3.2);

\draw (-4,0) -- (3,0);
\draw (0,0) -- (0,3);
\draw (0.447213595,-.2) -- (0.447213595,-.1);

\foreach \a/\b   in {
8.472135955/	8.972135955,
-2.618033989/	3.427050983,
1.618033989/	1.309016994,
0/	0.5,
0.618033989/	0.190983006,
0.381966011/	0.072949017,
0.472135955/	0.027864045,
0.437694101/	0.010643118,
0.450849719/	0.004065309
}
\draw[fill=gold!10]  (\a,\b) circle (\b)
;
\node at (2.8,2.7) {\Large$\displaystyle\frac{4\tau^3,\, 1}{2\tau^6,\, 8}$};
\node at (-2,	2) {\Large$\displaystyle\frac{-2\tau^2,\, 1}{~~2\tau^4,\, 2}$};
\node at (1.618033989, 1.309016994)  {\Large$\displaystyle\frac{2\tau,\,1}{2 \tau^2,\,2}$};
\node at (0,	0.45) {\large $\displaystyle\frac{0,\,1}{2,\,0}$};
\node at (0.618033989, 0.190983006) [scale=.87]{$\frac{2\varphi,1}{2\varphi^2,2}$};
\node at (2, .5) {\large$\dfrac{0,\, \hbox{\rm -} 1}{0,\, 0}$};
\draw [->, thick] (2.25,.5) to [bend left=45] (2.7,0.05) ;
\draw (-4,0) -- (3,0);
\draw [thick, ->] (0.447213595,-.5) -- (0.447213595,-0.05) node [below = 32pt, scale=.8] {$1/\sqrt{5}$};
\draw [thick, ->] (0,-.2) -- (0,-0.05) node [below = 11pt] {$0$};
\draw [thick, ->] (1,-.2) -- (1,-0.05) node [below = 11pt] {$1$};
\draw [thick, ->] (-2.618033989,-.2) -- (-2.618033989,-0.05) node [below = 11pt] {$-\varphi^2$};
\draw [thick, ->] (1.618033989,-.2) -- (1.618033989,-0.05) node [below = 11pt] {${}^{\phantom{2}}\varphi^{\phantom{2}}$};
\draw [thick, ->] (.618033989,-.2) -- (.618033989,-0.05) node [below = 11pt] {${}^{\phantom{2}}\tau^{\phantom{2}}$};
\draw [thick, ->] (.381966011,-.2) -- (.381966011,-0.05) node [below = 11pt] {$\tau^2~$};
\end{tikzpicture}
\end{center}
\caption{Half-plane Apollonian packing labeled with disk symbols.}
\label{fig:halfplane2}
\end{figure}
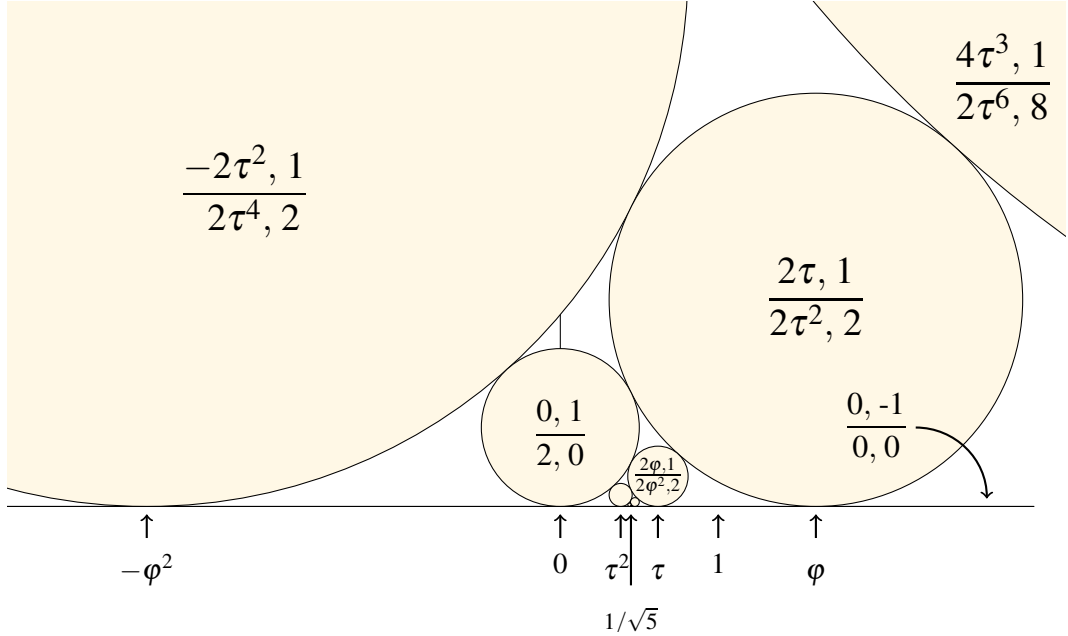

\begin{figure}
\begin{center}

\begin{tikzpicture}[scale=5.7]
\draw[fill=gold!50, draw=white] (-.55,-.05) rectangle (0.55, 1.05);
\draw[thick, fill=blue!20] (0,.5) circle (.5);
\draw (-.5,.5) -- (.5,.5);
\draw (0,0) -- (0,1);
\foreach \a/\b/\c   in {
0.000664308/	0.500740958/	0.000741838,
-0.001755194/	0.501962367/	0.00196848,
0.004600126/	0.505060138/	0.005101618,
-0.012865301/	0.514383844/	0.014696681,
0.033590288/	0.533590288/	0.035572686,
-0.177998211/	0.69900805/	0.233002683, 
0.259463815/	0.629731908/	0.209910636,
0.350372906/	0.324813547/	0.108271182,
0.295685999/	0.169413003/	0.056471001,
0.265216123/	0.234783877/	0.015652258,
0.253922357/	0.216106174/	0.006174462,
0.250918054/	0.223990141/	0.002262527
}
\draw[fill=gold!50, thick] (\a,\b) circle (\c)
;
\draw[fill=blue!20] (0,.25) circle (.25);

\draw [color=red]
(0.000664308,	0.500740958) -- 	
(-0.001755194,	0.501962367) --
(0.004600126,	0.505060138) --
(-0.012865301,	0.514383844) --
(0.033590288,	0.533590288) --
(-0.177998211,	0.69900805)	 --
(0.259463815,	0.629731908) -- (0.5, 0.75);
\draw [color=red]
(.5, .25) --
(0.350372906,	0.324813547) --
(0.295685999,	0.169413003) --
(0.265216123,	0.234783877) --
(0.253922357,	0.216106174) --
(0.250918054,	0.223990141); 
\end{tikzpicture}
%
\qquad  
\begin{tikzpicture}[scale=2.85]
\draw[fill=white, draw=white] (-1.1,-1.1) rectangle (0.6, 1.1);
\draw[fill=blue!20] (0,0) circle (1);
\draw (-1,0) -- (1,0);
\draw (0,-1) -- (0,1);
\draw [thick] (0,0) circle (1);
\foreach \a/\b/\c   in {
0.005267541/	0.994103732/	0.005882312,  
-0.013717028/	0.984520567/	0.01538388,
0.034669906/	0.960925249/	0.038449515,
-0.087352879/	0.895964206/	0.09978759,
0.186772685/	0.780159027/	0.197795447,
-0.355996422/	0.3980161/	0.466005367,  
0.51892763	/0.259463815/	0.419821272, 
0/	-0.5/	0.5,  
0.700745812/	-0.350372906/	0.216542365,
0.591371999/	-0.661173995/	0.112942002,
0.738496823/	-0.607745234/	0.043583863,
0.721753951/	-0.666542139/	0.017550414,
0.744347325/	-0.657699442/	0.006711776,
0.741887996/	-0.666663992/	0.002584,
0.74520866/	-0.665350052/	0.000987168
}
\draw [fill=gold!50] (\a,\b) circle (\c)    
;
\draw [color=red]
(0.005267541,	0.994103732)--
(-0.013717028,	0.984520567)--
(0.034669906,	0.960925249)--
(-0.087352879,	0.895964206)--
(0.186772685,	0.780159027)--
(-0.355996422,	0.3980161)--
(0.51892763,	0.259463815)--
(0,	-0.5)--
(0.700745812,	-0.350372906)--
(0.591371999,	-0.661173995)--
(0.738496823,	-0.607745234)--
(0.721753951,	-0.666542139)--
(0.744347325,	-0.657699442)--
(0.741887996,	-0.666663992)--
(0.74520866,	-0.665350052)
;
\end{tikzpicture}
\end{center}
\caption{Half-plane Apollonian packing under inversion -- two examples.}
\label{fig:wkole}
\end{figure}
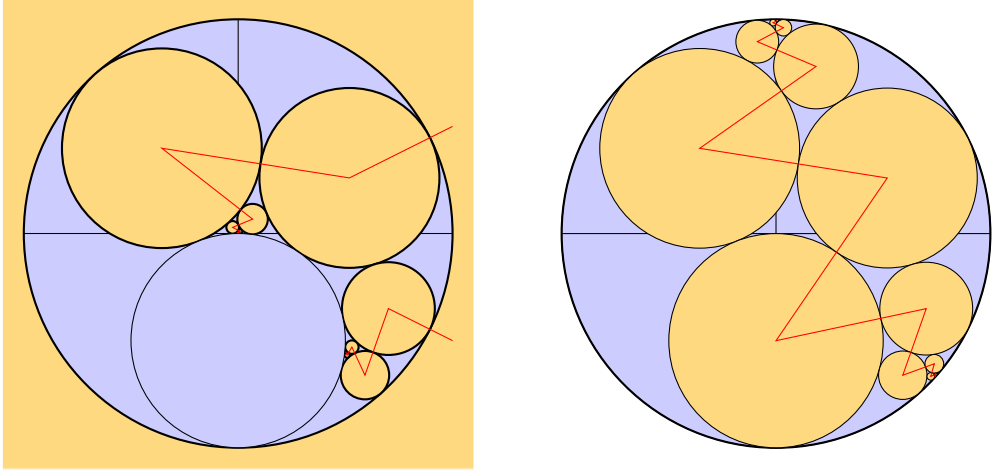

~\\
{\bf Remark 1:}
By inversion one may bring the half-plane arrangement \eqref{eq:HalfGold} into the interior of a disk.
Two such images are shown in Figure~\ref{fig:wkole}.  
The first (left), is obtained by inversion in a disk in the golden zigzag chain.
The image of the two limit disks for $n\to\infty$ and $n\to-\infty$ lie at the center $(0,0)$ 
and $(\frac{2\sqrt{5}}{9},\frac{5}{9})$, respectively. 
The second (right) arrangement can be obtained by inversion in a circle that is tangent to 
the horizontal axis line ($y=0$) at the point $x=5^{-1/2}$ and lies below it.  
The images of the limit disks are now at $(0,1)$ and at $(\frac{\sqrt{5}}{3},\frac{-2}{3})$, respectively.

\newpage

\section{Apollonian disk packing that fills the whole plane}

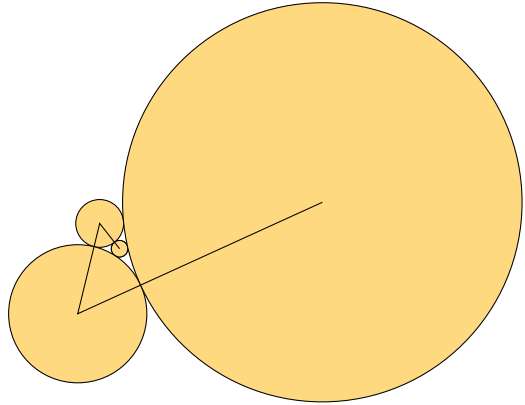
\begin{figure}[H]
\begin{center}
\begin{tikzpicture}[scale=.11]
\clip (-14, -19) rectangle  (51,32);   
\begin{scope}[rotate=0]  
\draw (-1,0) -- (1,0);
\draw (0,-1) -- (0,1);
\draw (0,0) circle (1);
\foreach \a/\b/\c   in {
0	/	0	/	1	,
-2.404185367	/	3.058171027	/	2.890053638	,
-5.058171027	/	-7.86654176	/	8.352410032	,
24.50341124	/	5.616741466	/	24.138913	
}	
\draw [fill=gold!50] (\a,\b) circle (\c)     
;
\draw 
(0,0	) --
(-2.404185367,	3.058171027	) --
(-5.058171027,	-7.86654176	) --
(24.50341124,	5.616741466	) 
;
\end{scope}
\end{tikzpicture}
\end{center}
\caption{The seed of the plane-filling arrangement}
\label{fig:plane-seed}
\end{figure}

Suppose four circles form a spiral (see Figure \ref{fig:plane-seed}) and their curvatures make a geometric sequence.
Then such an arrangement can be prolonged to infinity (in both directions, inward and outward) providing a skeleton for an 
Apollonian disk arrangement that fills the whole plane!

\noindent
Below, we inquire such a construction.
Denote the curvatures of the four consecutive disks in the spiral by $1$, $p$, $p^2$ and $p^3$.
Since they form a Descartes configuration, they must satisfy the Descartes formula~\eqref{eq:Descartes}:
\begin{equation}
(1+ p +  p^2 + p^3)^2 \ = \ 2(1^2 + p^2 + p^4 +p^6)\,.
\end{equation}
This is a sixth-degree polynomial equation
\begin{equation}
\label{eq:poly}
p^6  -2p^5 - p^4 - 4p^3 -p^2 -2p +1 = 0\,.
\end{equation}
Fortunately, it may be conveniently factorized:
\begin{equation}
(p^2 +1) \; (p^2 - 2\varphi p +1)\;(p^2 + 2\tau p +1) = 0 \, ,
\end{equation}
where $\varphi$ and $\tau$ are the golden ratio and golden cut, respectively.
Hence we end up with the following six roots:
\begin{equation}
 \pm i,  \quad -\tau \pm \sqrt{\tau}\,i, \quad  \varphi \pm\sqrt{\varphi}
\end{equation}
The only two real roots are mutual  reciprocals: 
$$
(\varphi +\sqrt{\varphi})\; (\varphi -\sqrt{\varphi}) = 1
$$
Thus, effectively, we obtain a unique solution and may choose either of the two real roots 
(the option is a choice between increasing or decreasing order in the chain).
Let us pick one and denote it as
\begin{equation}
\label{eq:rho}
\boxed{
\phantom{\int}
 \rho = \varphi + \sqrt{\varphi} 
\phantom{\int}}
\end{equation}
The radii in the spiral are the consecutive powers of $\rho\approx 2.89$.
Next, we need to find the angle through which the spiral turns at every vertex.
We get this result from the geometry of three consecutive mutually tangent disks 
(for instance $\rho^{-1} = \varphi-\sqrt{\varphi}$, \ $1$,  \ $\rho = \varphi+\sqrt{\varphi}$) 
by considering the triangle defined  by their centers. 
With the help of basic trigonometry, one finds the turn, expressed here by a unit complex number:
\begin{equation}
\label{eq:omega}
\boxed{
\phantom{\int}
\omega =  -\tau + \sqrt{\tau}\,i
\phantom{\int}}
\end{equation}
Surprisingly,  this ``turn number'' for the  the spiral is also a root of
the polynomial \eqref{eq:poly}!  
Note that $\omega$ is a unit complex number, $\omega\bar\omega = 1$; 
hence choosing the other, conjugated, root would only change the chirality of the disk spiral.   

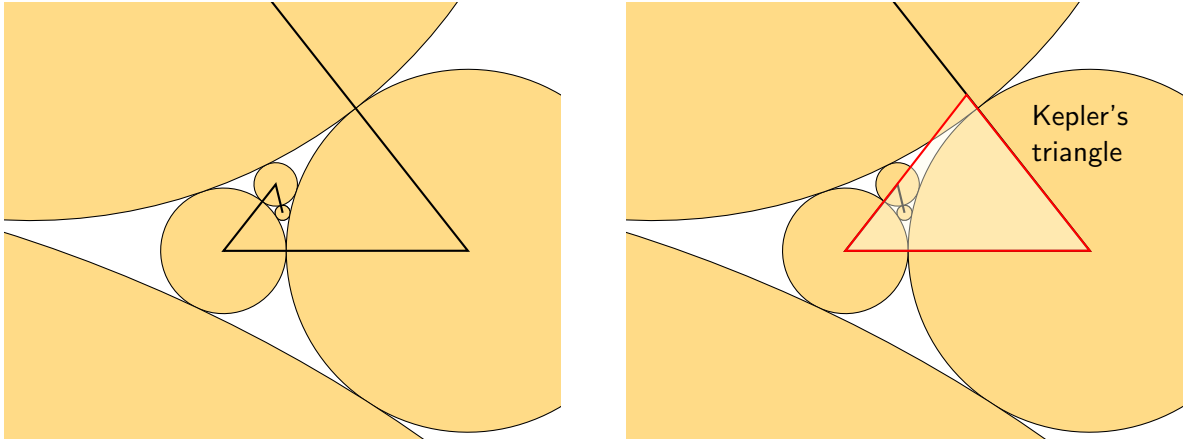
\begin{figure}
\begin{center}
\begin{tikzpicture}[scale=.1]

\clip (-37, -30) rectangle  (37,28);   
\begin{scope}[rotate=-24.5]
\draw (-1,0) -- (1,0);
\draw (0,-1) -- (0,1);
\draw (0,0) circle (1);
\foreach \a/\b/\c   in {
0	/	0	/	1	,
-2.404185367	/	3.058171027	/	2.890053638	,
-5.058171027	/	-7.86654176	/	8.352410032	,
24.50341124	/	5.616741466	/	24.138913	,
-58.93233408	/	48.69805233	/	69.76275334,
-7.785047432	/	-217.8193433	/	201.6180991
}	
\draw [fill=gold!47] (\a,\b) circle (\c)    
;
\draw [thick]
(0,0	) --
(-2.404185367,	3.058171027	) --
(-5.058171027,	-7.86654176	) --
(24.50341124,	5.616741466	) --
(-58.93233408,	48.69805233	) --
(-7.785047432,	-217.8193433	)
;
\end{scope}
\end{tikzpicture}
\qquad
\begin{tikzpicture}[scale=.1]
\clip (-37, -30) rectangle  (37,28);   
\begin{scope}[rotate=-24.5]
\draw (-1,0) -- (1,0);
\draw (0,-1) -- (0,1);
\draw (0,0) circle (1);
\foreach \a/\b/\c   in {
0	/	0	/	1	,
-2.404185367	/	3.058171027	/	2.890053638	,
-5.058171027	/	-7.86654176	/	8.352410032	,
24.50341124	/	5.616741466	/	24.138913	,
-58.93233408	/	48.69805233	/	69.76275334,
-7.785047432	/	-217.8193433	/	201.6180991
}	
\draw [fill=gold!47] (\a,\b) circle (\c)
;
\draw [thick]
(0,0	) --
(-2.404185367,	3.058171027	) --
(-5.058171027,	-7.86654176	) --
(24.50341124,	5.616741466	) --
(-58.93233408,	48.69805233	) --
(-7.785047432,	-217.8193433	)
;
\draw [color=red, thick, fill= yellow!10, fill opacity = .42]
(-5.058171027,	-7.86654176	) --
(24.50341124,	5.616741466	) --
(1, 17.7) -- cycle
;
\node at (17, 19	) [align=left]  {\sf Kepler's \\  \sf triangle};
\end{scope}
\end{tikzpicture}
\end{center}
\caption{Left: disks filling the whole plane. Right: Kepler's triangles.}
\label{fig:plane}
\end{figure}

~
The findings \eqref{eq:rho} and \eqref{eq:omega} 
imply that that the positions of the centers of the disks, understood as complex numbers , 
satisfy the recurrence:
$$
z_{n+1} = z_{n} + (1+\rho )\, (\rho\omega)^n 
$$
Choosing $z_0=0$,  we obtain  $z_n=\sum_{k=0}^n z_k$, which 
-- as a finite geometric series -- may easily be algebraically reduced.
Here is the result:    


\newpage
\noindent
{\bf Theorem 2:}  Let $D_n$ be a chain of disks in the complex plane $\mathbb C\cong \mathbb R^2$, 
described below:

\begin{equation}
\label{eq:spiral}
\begin{array}{rll}
\hbox{radius:} & r_n = \rho^n     &\in\mathbb R\\
\hbox{center:} &z_{n} = (1+\rho )\, \displaystyle{\frac{(\rho\omega)^n -1}{\rho\omega -1}}\quad &\in \mathbb C
\end{array}
\end{equation}
where:
$$
\rho = \varphi +\sqrt{\varphi}   \,,\qquad
\omega = -\tau + \sqrt{\tau}\, i  \,,\qquad \hbox{and}\qquad
\rho\,\omega = (1+\sqrt{\tau})(-1+ \sqrt{\varphi}\,i)
$$
%
%
Then any quadruple of consecutive disks $(D_n,  D_{n+1}, D_{n+2},  D_{n+3})$
forms  a Descartes configuration.
The disks form a spiral-like arrangement
that may be completed to form an unbounded Apollonian disk packing that fills the whole plane
with disks of positive curvatures. 

~\\
{\bf Proof:} Although the derivation preceding the theorem validates this statement, 
one may still verify the Descartes formula for the entries, a monotonous but not difficult task.
The unboundedness is evident.
\QED

~\\
{\bf Corollary 3:} The center of the spiral lies at
$$
z_{-\infty} \ = \ \lim_{n\to\-\infty} z_n \  = \  \frac{1+\varrho}{1-\varrho\omega} \ \approx \ 0.84+0.68 i
\ \in \mathbb C
$$
If the spiral's center is shifted to the coordinate origin,
the equations defining the spiral become
\begin{equation}
\begin{array}{rl}
\hbox{radius:} & r_n = \rho^n\\
\hbox{center:} &z_{n} = (1+\rho )\, \displaystyle{\frac{(\rho\omega)^n}{\rho\omega -1}}\\
\end{array}
\end{equation}

 ~\\
{\bf Corollary 4:} The symbols for the seed of this construction may now be found:
$$
\left( ~
\frac{1+\bar \rho,\ 0}{\bar\rho,\ \bar\rho +2}, 
\quad
\frac{0,\ 0}{1,\ -1}, 
\quad
\frac{-(1+\rho)\bar \omega ,\ 0}{\bar\rho,\ \rho +2}
~ \right)\,,
$$
where $\bar \varrho = \varphi - \sqrt{\varphi}$.
 In terms of real numbers and golden means,
the seed is
\begin{equation}\left( ~
\dfrac{\!\varphi\!-\!\sqrt{\varphi}\!+\!1, \ \ 0}
       {\varphi\!-\!\sqrt{\varphi},\ \varphi\!-\!\sqrt{\varphi}\!+\!2} ,
\quad
\dfrac{0,\ \  0}{1,\ -1}
\quad
\dfrac{\tau\! +\!\sqrt{\tau}\!+\!1, \ \sqrt{\varphi}+\sqrt{\tau}+1}
       {\varphi\! +\!\sqrt{\varphi},\ \varphi\!+\!\sqrt{\varphi}\!+\!2}
       ~ \right)\,.
\end{equation}

~\\
{\bf Remark 2:}
The angle formed by two oriented arms joining the centers of three consecutive disks 
make an angle $\theta$ that satisfies 
$$
-\bar\omega =  \tau + \sqrt{\tau}\,i = e^{\theta i} \qquad \Rightarrow \qquad 
\cos \theta = -\tau \quad \hbox{and} \quad  \tan\theta = \sqrt{\varphi} 
$$ 
This happens to be the angle that lies at the base of the Kepler right triangle 
$(1, \sqrt{\varphi}, \varphi)$, which Johannes Kepler found as the right triangle
with the sides that form a geometric progression.
The silhouette of the Khufu Pyramid is effectively made of two such triangles juxtaposed
(intentionally or not).
Hence three segments joining the centers of any three consecutive circles of the chain (\eqref{eq:spiral} form, 
under closing, such a triangle (see Figure \ref{fig:plane} right).

~

In words:  the ``golden spiral'' is a chain consisting of disks, the ratio of the radii of any two consecutive being
$$
\rho = \varphi+\sqrt{\varphi} \approx 2.89005
$$
They are arranged into a spiral. 
The angle made by arms drawn from any disk to two neighboring disks in the chain is 
$$
\theta = \arccos \tau = \arctan \sqrt{\varphi} \approx 51.83^\circ   \,.
$$

\section{Concluding remarks}

The Apollonian gasket, understood as an arrangement of {\it circles} may or may not be bounded in the plane.  
But as an arrangement of disks, {\it every} Apollonian gasket fills the entire plane.  
This is because a disk with negative curvature is unbounded and fills the region outside the circle.
Denote by $\mathcal C(A)$ the set of curvatures of an Apollonian packing $A$.
Also, denote by $A_0\subset A$ the subset of disks of curvature 0.
Here is a list of the types of Apollonian disk packings
that accounts for four visually different forms:

\begin{figure}[H]
\centering
\begin{tikzpicture}[scale=.9]

\draw (0,0) circle (1);

\foreach \a/\b/\c   in {
1 / 0 / 2 
}
\draw (\a/\c,\b/\c) circle (1/\c)
          (-\a/\c,\b/\c) circle (1/\c);
\foreach \a/\b/\c in {
0 / 2 / 3 ,
0 /4 /15 ,
0 / 6 / 35 , 
0 / 8/ 63
}
\draw (\a/\c,\b/\c) circle (1/\c)
          (\a/\c,-\b/\c) circle (1/\c) ;
\foreach \a/\b/\c/\d in {
3 / 4 /6, 	8 / 6 / 11,	5 / 12/ 14,	15/ 8 / 18,	8 / 12 / 23,	7 / 24 / 26,
24/	10/	27, 	21/	20/	30, 	16/	30/	35, 	3/	12/	38, 	35/	12/	38, 	24/	20/	39, 	9/	40/	42,
16/	36/	47, 	15/	24/	50, 	 48/	14/	51, 	45/	28/	54, 	24/	30/	59, 	40/	42/	59, 	11/	60/	62
}
draw (\a/\c,\b/\c) circle (1/\c)       (-\a/\c,\b/\c) circle (1/\c)
          (\a/\c,-\b/\c) circle (1/\c)       (-\a/\c,-\b/\c) circle (1/\c) 
;
\node at (-1/2,0) [scale=1.7, color=black] {\sf 2};
\node at (1/2,0) [scale=1.7, color=black] {\sf 2};
\node at (0,2/3) [scale=1.4, color=black] {\sf 3};
\node at (0,-2/3) [scale=1.4, color=black] {\sf 3};
\end{tikzpicture}
                                                         \qquad
%
\begin{tikzpicture}[scale=.67, rotate=90, shift={(0,2cm)}]  
\clip (-1.25,-2.1) rectangle (1.1,2.2);
\draw (1,-2) -- (1,3);
\draw (-1,-2) -- (-1,3);
\draw (0,0) circle (1);
\draw (0,2) circle (1);
\draw (0,-2) circle (1);
\foreach \a/\b/\c in {
3/4/4,  5/12/12,  7/24/24, 9/40/40  
}
\draw (\a/\c,\b/\c) circle (1/\c)    (\a/\c,-\b/\c) circle (1/\c)
          (-\a/\c,\b/\c) circle (1/\c)    (-\a/\c,-\b/\c) circle (1/\c)   ;
\foreach \a/\b/\c in {
8/    6/   9, 	
15/   8/   16 , 	
24/  20/  25, 	
24/  10/  25, 	
21/  20/  28,
16/	30/	33,    
35/	12/	36,
48/	42/	49,
48/	28/	49,
48/	14/	49
}
\draw (\a/\c, \b/\c) circle (1/\c)          (-\a/\c, \b/\c) circle (1/\c)
          (\a/\c,-\b/\c) circle (1/\c)         (-\a/\c,-\b/\c) circle (1/\c)
          (\a/\c,2-\b/\c) circle (1/\c)       (-\a/\c,2-\b/\c) circle (1/\c)
          (\a/\c,2+\b/\c) circle (1/\c)       (-\a/\c,2+\b/\c) circle (1/\c)
          (\a/\c,-2+\b/\c) circle (1/\c)       (-\a/\c,-2+\b/\c) circle (1/\c)
;
\node at (0,0) [scale=1.9, color=black] {\sf 1};
\node at (0,-2) [scale=1.9, color=black] {\sf 1};
\node at (0,2) [scale=1.9, color=black] {\sf 1};
\node at (-3/4,-1) [scale=1.1, color=black] {\sf 4};
\end{tikzpicture}
                                                       \qquad
\begin{tikzpicture}[scale=.6]

\clip (-3.4, -.05) rectangle  (3.2,2.8);
\draw (-4,0) -- (4,0);
\foreach \a/\b   in {
8.472135955/	8.972135955,
-2.618033989/	3.427050983,
1.618033989/	1.309016994,
0/	0.5,
0.618033989/	0.190983006,
0.381966011/	0.072949017,
0.472135955/	0.027864045,
0.437694101/	0.010643118,
0.450849719/	0.004065309
}
\draw[fill=gold!10]  (\a,\b) circle (\b)
;
\node at (-2.618033989+1,	3.427050983-2) [align=left,scale=1.6] {$\varphi^{-2}$ }; 
\node at (1.618033989, 1.309016994)  [align=left, scale=1.6] {$1$ }; 
\node at (0,	0.5) [align=left, scale=.9] {$\varphi^{2}$ }; 
\end{tikzpicture}
                                                      \qquad
\begin{tikzpicture}[scale=.032]
\clip (-47, -30) rectangle  (37,28);   
\begin{scope}[rotate=-24.5]
\draw (-1,0) -- (1,0);
\draw (0,-1) -- (0,1);
\draw (0,0) circle (1);
\foreach \a/\b/\c   in {
0	/	0	/	1	,
-2.404185367	/	3.058171027	/	2.890053638	,
-5.058171027	/	-7.86654176	/	8.352410032	,
24.50341124	/	5.616741466	/	24.138913	,
-58.93233408	/	48.69805233	/	69.76275334,
-7.785047432	/	-217.8193433	/	201.6180991
}	
\draw [fill=gold!10] (\a,\b) circle (\c)
;
\node at (-10,	-30	) [scale=.8] {$\varphi\!+\!\sqrt{\varphi}$};
\node at (24.50341124,	5.616741466	) [scale=1.4]{$1$};
\node at (-25,	5	) {~~~$\varphi\!-\!\sqrt{\varphi}$};
\draw [->] (-18,-27) to [bend left=45] (-6,-9);
\end{scope}
\end{tikzpicture}
\vspace {-.1in}
$$
\hbox{Type A} \qquad\qquad\qquad
\hbox{Type B} \qquad\qquad\qquad\qquad
\hbox{Type C} \qquad\qquad\qquad\qquad\quad\quad
\hbox{Type D}
$$
\vspace {-.3in}
\caption{The four types of Apollonian disc packing.}
\label{fig:types}
\end{figure}
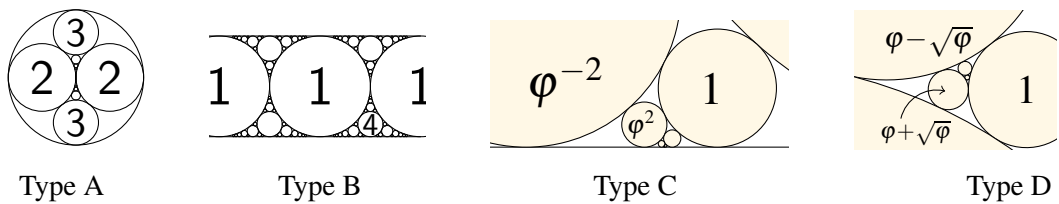

\begin{itemize}
\item[]
{\bf Type A}:   $\inf \mathcal C(A)=\min \mathcal C(A) <0$.
There exists a circle with negative curvature.  This is the most frequently discussed case.  

\item[]
{\bf Type B}: $\inf \mathcal C(A) =\min \mathcal C(A) = 0$ and $\hbox{\rm card}\, A_0= 2$.
There exist exactly two 0-curvature disk.  
This is the example of an Apollonian Belt, related to the arrangement known as the Ford circles. 
In this case, 0 is not an accumulation point of $\mathcal C(A)$. 
\item[]
{\bf Type C}:  $\inf \mathcal C(A) = 0$ and $\hbox{\rm card}\, A_0= 1$.
There exists exactly one 0-curvature disk.  
In this case 0 is an accumulation point of $\mathcal C(A)$ and $\inf\mathcal C(A) \in \mathcal C(A)$.
This is a case of the half-plane filling.  
This case was also discussed in \cite{CD}.  
\item[]
{\bf Type D}: $\inf \mathcal C(A) = 0$ and $\hbox{\rm card}\,  A_0= 0$.
There are no 0-curvature disks but  0 is an accumulation point of $\mathcal C(A)$,
$\inf\mathcal C(A) \not\in \mathcal C(A)$.
This is the case of a plane filing arrangement.
\end{itemize}

The examples of type {\bf A} and {\bf B} are quite popular.
The unbounded arrangements {\bf C} and {\bf D} are similar and 
the main difference between them is whether  $\mathcal C(A)$ contains its $\inf \mathcal C(A)$ .
\\

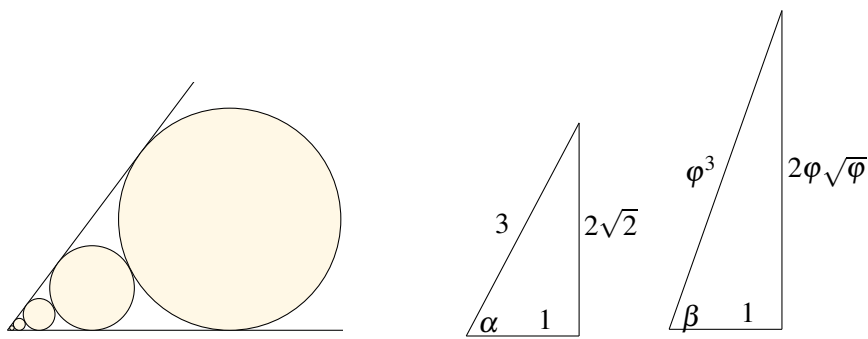
\begin{figure}[H]
\centering
\begin{tikzpicture}[scale=1.5] 
\clip (-.25,-.1) rectangle (2.97,2.2);
\draw (0,0) -- (3,0);
\draw (0,0) -- (3,4);
\draw (3,0) -- (3,4);
\draw [fill=gold!10] (2-.03,1-.015) circle (1-.015);
\draw [fill=gold!10] (.75,.75/2) circle  (.75/2);
\draw [fill=gold!10] (.28115,.28125/2) circle (.28125/2);
\draw [fill=gold!10] (.1055,.1055/2) circle (.1055/2);
\draw [fill=gold!10] (.04,.04/2) circle (.04/2);
\end{tikzpicture}
\qquad\qquad
\begin{tikzpicture}[scale=1.5, rotate=0, shift={(0,0cm)}]  
\draw (0,0) -- (1,0);
\draw (1,0) -- (1,1.89);
\draw (0,0) -- (1,1.89);
\node at (1.28,1) [scale=1, color=black] {$2\sqrt{2}$};
\node at (0.32,1) [scale=1, color=black] {$3$};
\node at (0.7,.15) [scale=1, color=black] {$1$};
\node at (0.2,.1) [scale=1, color=black] {$\alpha$};
\end{tikzpicture}
\ \
\begin{tikzpicture}[scale=1.5, rotate=0, shift={(0,0cm)}]  
\draw (0,0) -- (1,0);
\draw (1,0) -- (1,2.828);
\draw (0,0) -- (1,2.828);
\node at (1.4,1.4) [scale=1, color=black] {$2\varphi\sqrt{\varphi}$};
\node at (0.28,1.4) [scale=1, color=black] {$\varphi^3$};
\node at (0.7,.15) [scale=1, color=black] {$1$};
\node at (0.2,.1) [scale=1, color=black] {$\beta$};
\end{tikzpicture}                                     
%

\caption{Left: An example of a wedge chain of disks in a triangle;  
Center: the triangle for chain B; 
Right: The triangle for chain B.}
\label{fig:3triangles}
\end{figure}

Note that chains of disks with radii that form a geometric progression may be inscribed (and determined by) an angle,
as in Figure~\ref{fig:3triangles}, left.  The angle may be conveniently described by a right triangle. 

%
%
%
%
%

~\\\\
{\bf Corollary 5:} Each of the golden chains \eqref{eq:HalfGold} and \eqref{eq:spiral} can be inscribed 
in the base angle $\alpha$ and $\beta$, respectively, of the triangle of proportions 
shown in Figure \ref{fig:3triangles}.
In particular:
\begin{equation}
\begin{array}{ccc}
 \cos\alpha \ = \ \dfrac{1}{3}\,,    &\qquad &  \cos\beta \ = \ \dfrac{1}{\varphi^3} \\[7pt]
\hbox{\footnotesize \sf [zigzag]}&&\hbox{\footnotesize \sf [spiral]}
\end{array}  
\end{equation}

~\\
{\bf Proof:} Use trigonometric identities for of doubling an angle. The first result is straightforward; 
the second, despite its simple appearance, requires somewhat tedious algebraic manipulations
that are left to the reader as an exercise.
\QED

~\\

%
%
%


~\\



\begin{thebibliography}{9}


\bibitem{CD}
Michael Ching and John R. Doyle,
Apollonian circle packing of the half-plane,
arXiv:1102.1628v2 [math.MG]


\bibitem
{jk1} Jerzy Kocik, 
        A theorem on circle configurations. arXiv:0706.0372v2. 







\bibitem{jk}
J. Kocik,   
``On a Diophantine equation that generates all integral Apollonian gaskets,''
ISRN Geometry, 348618 (2012)

\bibitem
{jk3}  Jerzy Kocik, 
Pages on Apollonian packings,     
{\sf http://Lagrange.math.siu.edu/Kocik/apollo/apollo.html}

\bibitem
{LMW} 
Jeffrey C. Lagarias, Colin L. Mallows, Allan R. Wilks,
        Beyond the Descartes circle theorem. {\it Am. Math. Monthly} {\bf 109} (2002), 338--361.


\bibitem{N} 
S. Northshield, 
``On integral Apollonian circle packings,'' 
J. Number Theory,
{\bf 119}(2), 171-193 (2006).



\end{thebibliography}
\end{document}